\documentclass{amsart}
\usepackage{pictex}
\usepackage{amsmath}
\usepackage{txfonts}
\theoremstyle{definition}

\theoremstyle{remark}

\numberwithin{equation}{section}

\begin{document}

\title{Geometrical Equivalents of Goldbach Conjecture\\
 and \\
Fermat Like Theorem}

% Remove or comment out any unused author tags.
% author one information
\author{K. K. Nambiar}
\address{Formerly, Jawaharlal Nehru University, New Delhi, 110067, India}
\curraddr{1812 Rockybranch Pass, Marietta, Georgia, 30066-8015}
\email{knambiar@fuse.net}
%\thanks{}

% author two information
%\author{}
%\address{}
%\curraddr{}
%\email{}
%\thanks{}

% Use this \subjclass if you are using amsart version 2.0 (December 1999).
% Use this one if you are using an older version of amsart.
%\subjclass{}
\date{November 02, 2002}

% at present the "communicated by" line appears only in ERA and PROC
\commby{}

%\dedicatory{Dedicated to Professor Stephen Hawking who lit up 
%the black hole for us.}

\begin{abstract}
Five geometrical equivalents of Goldbach conjecture are 
given, calling one of them Fermat Like Theorem.
\end{abstract}

\maketitle
%\tableofcontents
%\end{document}

%-----------------------------------------------------------------------
% End of journal.top
%-----------------------------------------------------------------------

\section {Introduction}

\noindent 
The well-known Goldbach conjecture \cite {BalCox:MRE} 
states that every even 
number, greater than $4$, can be written as the sum of two odd 
primes. The purpose of this paper is to state a few geometrical 
equivalents of Goldbach conjecture. 

\section {Goldbach Circle} 

All the different forms of the conjecture can be obtained 
by considering the \emph {Goldbach Circle} as given in 
Figure 1. Here, $ADB$ is a right-angled triangle, inscribed in 
a circle of radius $n$, an integer greater than $2$. 

%\begin{figure}
\vglue 8truept 
\beginpicture
\setplotarea x from 0 to 300, y from 0 to 110

\linethickness = 1pt 
\setlinear
\plot 135 65 235 65 /
\plot 171 65 171 113 / 
\plot 135 65 171 113 /
\plot 171 113 235 65 /
\plot 185 65 171 113 /

\put {$A$} at 130 65
\put {$B$} at 240 65
\put {$C$} at 185 60
\put {$D$} at 171 118
\put {$E$} at 171 60

\circulararc 360 degrees from 185 115 center at 185 65 

\put {Figure~1.\enskip Goldbach Circle} at 185 5
\endpicture
\vskip 8truept
%\end{figure}

Further, diameter $AB=2n$, radius $CD=n$, $AE=p_1$ a prime, $EB=p_2$ 
a prime, $DE=\sqrt {p_{1}p_{2}}$, $AD=\sqrt{p_{1}(p_{1}+p_{2})}$, 
$BD=\sqrt{p_{2}(p_{1}+p_{2})}$, and $EC={(p_{2}-p_{1})}/2$.
From the figure, we can write several different versions of 
Goldbach conjecture.

\begin{description}
\item [Version 1] For every $n$, there is a right-angled triangle 
$ADE$ with sides, \\
$AD=\sqrt {2np_{1}}$, $AE=p_{1}$, and 
$DE=\sqrt{p_{1}p_{2}}$. 
\item [Version 2] For every $n$, there is a right-angled triangle 
$BDE$ with sides, \\
$BD=\sqrt {2np_{2}}$, $BE=p_{2}$, and 
$DE=\sqrt{p_{1}p_{2}}$. 
\item[Version 3] For every $n$, there is a right-angled triangle 
$DEC$ with sides, \\
$CD=n$, $DE=\sqrt{p_{1}p_{2}}$, and 
$EC={(p_{2}-p_{1})}/2$. 
\item [Version 4] For every $n$, there is a right-angled triangle
$ADB$ with sides, \\
$AB=2n$, $AD=\sqrt {2np_{1}}$, and 
$BD=\sqrt {2np_{2}}$.
\item [Version 5] Corresponding to every $n$, there is a Goldbach 
Circle.
\end{description}

\section{Conclusion}
If we call $g=\sqrt {p_{1}p_{2}}$ a \emph {geometric prime} 
and $h={(p_{2}-p_{1})}/2$ a \emph {difference prime}, we can 
write version 3 as
\[
n^{2}= g^{2}+h^{2},
\]
which looks similar to Fermat's last theorem in form. Hence, we 
may call this version \emph {Fermat Like Theorem} which, 
of course, is only a conjecture.

%\bibliography{references}
\bibliographystyle{amsplain}

\end{document}